\documentclass[12pt]{article}                                               

\input{amssym.def}                                                           
\input{amssym.tex}                                                           
                                                                             
\begin{document}                                                             
\title{Toric $AF$-algebras  and faithful representation
of the mapping class groups}

\author{Igor Nikolaev
\footnote{Partially supported 
by NSERC.}}

\date{}
 \maketitle

\newtheorem{thm}{Theorem}
\newtheorem{lem}{Lemma}
\newtheorem{dfn}{Definition}
\newtheorem{rmk}{Remark}
\newtheorem{cor}{Corollary}
\newtheorem{prp}{Proposition}
\newtheorem{exm}{Example}
\begin{abstract}
There exists a covariant non-injective functor
from the space of generic Riemann surfaces to the so-called
toric $AF$-algebras;  such a functor maps isomorphic  Riemann surfaces 
to the  stably isomorphic toric  $AF$-algebras.  We use the functor  to 
construct a faithful representation of the mapping class group of surface of genus 
$g>1$  into the matrix group $GL_{6g-6}({\Bbb Z})$.

\vspace{7mm}


\vspace{5mm}
{\it AMS (MOS) Subj. Class.:  14H55, 20F65, 46L85}
\end{abstract}

\section{Introduction}
{\bf A. The mapping class group.}
The mapping class group has  been introduced in the 1920-ies by M.~Dehn \cite{Deh1}. 
Such a  group, $Mod~(X)$, is defined as the group of isotopy classes of the
orientation-preserving diffeomorphisms of a two-sided closed surface $X$ of
genus $g\ge 1$.  The group is known to be  prominent in  algebraic geometry \cite{HaLo1}, 
topology \cite{Thu1}  and dynamics \cite{Thu2}.  When $X$
is a   torus, the $Mod~(X)$ is isomorphic to the group  $SL_2({\Bbb Z})$. 
(The $SL_2({\Bbb Z})$ is called a modular group, hence our notation for the mapping class groups.) 
A little is known about the representations of $Mod~(X)$ beyond the case $g=1$.
Recall, that the group is called {\it linear},  if there exists a faithful
representation into the matrix group $GL_m(R)$, where $R$ is a commutative ring.
The braid groups are known  to be linear \cite{Big1}.  Using a modification of the argument for 
the braid groups, it is possible to prove, that $Mod~(X)$ is linear in the case $g=2$
\cite{BiBu1}. Whether the mapping class group is linear for  $g\ge 3$, 
is an open problem, known as a {\it Harvey conjecture} \cite{Har1}, p.267.

\medskip\noindent
{\bf B. The toric $AF$-algebras.}
 Denote by $T_S(g)$ the Teichm\"uller
space of genus $g\ge 1$ with a distinguished point $S$. Let
$q\in H^0(S, \Omega^{\otimes 2})$ be a holomorphic quadratic 
differential on the Riemann surface $S$, such that all zeroes
of $q$ (if any) are simple. By $\widetilde S$ we mean a double 
cover of $S$ ramified over the zeroes of $q$ and by
$H_1^{odd}(\widetilde S)$ the odd part of the integral homology of $\widetilde S$ relative to the  zeroes.
Note that $H_1^{odd}(\widetilde S)\cong {\Bbb Z}^{n}$, where $n=6g-6$ if $g\ge 2$ and $n=2$ if $g=1$. 
The fundamental result of Hubbard and Masur \cite{HuMa1} implies, that
$T_S(g)\cong Hom~(H_1^{odd}(\widetilde S); {\Bbb R})-\{0\}$, where $0$ is the zero homomorphism. 
Finally, denote by $\lambda=(\lambda_1,\dots,\lambda_{n})$ the image of a basis of 
$H_1^{odd}(\widetilde S)$ in the real line ${\Bbb R}$, such that $\lambda_1\ne0$. 
Note that such an option always exists, since the zero homomorphism is excluded.
We let $\theta=(\theta_1,\dots,\theta_{n-1})$, where $\theta_i=\lambda_{i-1}/\lambda_1$. 
Recall that,  up to a scalar multiple, the vector $(1,\theta)\in {\Bbb R}^{n}$ is the limit
of a generically convergent Jacobi-Perron continued fraction \cite{B}:
$$
\left(\matrix{1\cr \theta}\right)=
\lim_{k\to\infty} \left(\matrix{0 & 1\cr I & b_1}\right)\dots
\left(\matrix{0 & 1\cr I & b_k}\right)
\left(\matrix{0\cr {\Bbb I}}\right),
$$
where $b_i=(b^{(i)}_1,\dots, b^{(i)}_{n-1})^T$ is a vector of the non-negative integers,  
$I$ the unit matrix and ${\Bbb I}=(0,\dots, 0, 1)^T$. We introduce an  $AF$-algebra, 
${\Bbb A}_{\theta}$, via the Bratteli diagram \cite{E}, shown in Fig.1.
(The numbers $b_j^{(i)}$ of the diagram indicate the multiplicity  of edges of 
the graph.)  Let us call ${\Bbb A}_{\theta}$ a {\it toric $AF$-algebra}.

\begin{figure}[here]
\begin{picture}(300,200)(0,0)

\put(40,110){\circle{3}}

\put(40,110){\line(1,1){60}}
\put(40,110){\line(2,1){60}}
\put(40,110){\line(1,0){60}}
\put(40,110){\line(2,-1){60}}
\put(40,110){\line(1,-1){60}}
\put(40,110){\line(2,-3){60}}

\put(100,20){\circle{3}}
\put(100,50){\circle{3}}
\put(100,80){\circle{3}}
\put(100,110){\circle{3}}
\put(100,140){\circle{3}}
\put(100,170){\circle{3}}

\put(160,20){\circle{3}}
\put(160,50){\circle{3}}
\put(160,80){\circle{3}}
\put(160,110){\circle{3}}
\put(160,140){\circle{3}}
\put(160,170){\circle{3}}

\put(220,20){\circle{3}}
\put(220,50){\circle{3}}
\put(220,80){\circle{3}}
\put(220,110){\circle{3}}
\put(220,140){\circle{3}}
\put(220,170){\circle{3}}

\put(160,20){\line(2,1){60}}
\put(160,19){\line(1,0){60}}
\put(160,21){\line(1,0){60}}
\put(160,50){\line(2,1){60}}
\put(160,49){\line(2,-1){60}}
\put(160,51){\line(2,-1){60}}
\put(160,80){\line(2,1){60}}
\put(160,79){\line(1,-1){60}}
\put(160,81){\line(1,-1){60}}
\put(160,110){\line(2,1){60}}
\put(160,109){\line(2,-3){60}}
\put(160,111){\line(2,-3){60}}
\put(160,140){\line(2,1){60}}
\put(160,139){\line(1,-2){60}}
\put(160,141){\line(1,-2){60}}
\put(160,170){\line(2,-5){60}}


\put(100,20){\line(2,1){60}}
\put(100,19){\line(1,0){60}}
\put(100,20){\line(1,0){60}}
\put(100,21){\line(1,0){60}}
\put(100,50){\line(2,1){60}}
\put(100,49){\line(2,-1){60}}
\put(100,51){\line(2,-1){60}}
\put(100,80){\line(2,1){60}}
\put(100,79){\line(1,-1){60}}
\put(100,81){\line(1,-1){60}}
\put(100,110){\line(2,1){60}}
\put(100,109){\line(2,-3){60}}
\put(100,111){\line(2,-3){60}}
\put(100,140){\line(2,1){60}}
\put(100,139){\line(1,-2){60}}
\put(100,141){\line(1,-2){60}}
\put(100,170){\line(2,-5){60}}


\put(250,20){$\dots$}
\put(250,50){$\dots$}
\put(250,80){$\dots$}
\put(250,110){$\dots$}
\put(250,140){$\dots$}
\put(250,170){$\dots$}


\put(125,5){$b_5^{(1)}$}
\put(100,30){$b_4^{(1)}$}
\put(93,60){$b_3^{(1)}$}
\put(90,90){$b_2^{(1)}$}
\put(88,120){$b_1^{(1)}$}

\put(185,5){$b_5^{(2)}$}
\put(160,30){$b_4^{(2)}$}
\put(153,60){$b_3^{(2)}$}
\put(150,90){$b_2^{(2)}$}
\put(148,120){$b_1^{(2)}$}


\end{picture}

\caption{The toric $AF$-algebra ${\Bbb A}_{\theta}$ of genus $g=2$.}
\end{figure}
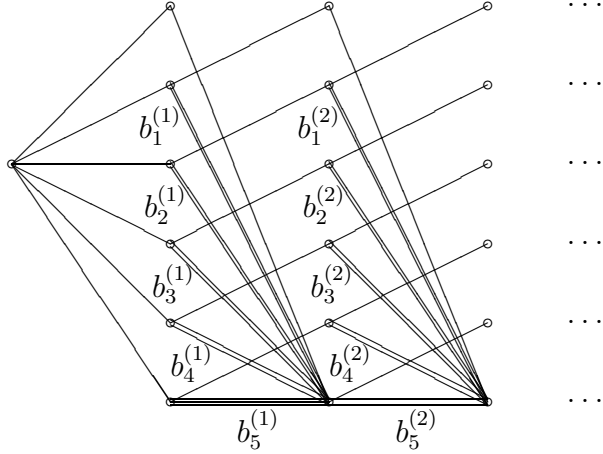

\medskip\noindent
{\bf C. The Teichm\"uller functor.}
Let $R,R'\in T_S(g)$ be a pair of isomorphic  Riemann surfaces;  let 
${\Bbb A}_{\theta}$ and ${\Bbb A}_{\theta'}$ be the corresponding  toric $AF$-algebras.  
We look for an answer to the following  elementary question: {\it How  are the algebras ${\Bbb A}_{\theta}$ and  
${\Bbb A}_{\theta'}$ related to each other?}    Recall,  that the stable isomorphism between the $C^*$-algebras is a
fundamental equivalence in noncommutative geometry; the $C^*$-algebras
${\Bbb A}$ and ${\Bbb A}'$ are said to be {\it stably isomorphic}, whenever ${\Bbb A}\otimes {\cal K}$
is isomorphic to  ${\Bbb A}'\otimes {\cal K}$, where ${\cal K}$ is the $C^*$-algebra
of compact operators.  
Denote by $V$ the maximal subset of $T_S(g)$,  such that for each Riemann surface $R\in V$, 
there exists a convergent Jacobi-Perron continued fraction. 
Let $F$ be the map which sends the Riemann surfaces into the toric $AF$-algebras
according to the formula  $R\mapsto {\Bbb A}_{\theta}$;
we shall call $F$ the {\it Teichm\"uller functor}. 
 Let $W$ be the image of  $V$ under $F$;
the following lemma relates the algebras ${\Bbb A}_{\theta}$ and ${\Bbb A}_{\theta'}$.
\begin{lem}\label{lm0}
{\bf (\cite{Nik1})}
The set $V$ is a generic subset of $T_S(g)$ and the  map $F$ has the following properties:
(i) $V\cong W\times (0,\infty)$ is a trivial fiber bundle, whose
projection  map $p: V\to W$ coincides with $F$;
(ii) $F: V\to W$ is a covariant functor, which maps isomorphic Riemann surfaces 
$R,R'\in V$ to  stably isomorphic toric $AF$-algebras ${\Bbb A}_{\theta},{\Bbb A}_{\theta'}\in W$.
\end{lem}

\medskip\noindent
{\bf D. The result.}
Recall  that  $Mod~(X)$ acts on $T_S(g)$ by 
isomorphisms of the Riemann surfaces; the action is properly
discontinuous and free for a  finite index subgroup  of $Mod~(X)$ \cite{HaLo1}.
Lemma \ref{lm0} extends the action to the toric $AF$-algebras,
where $Mod~(X)$  acts by  the stable isomorphisms;  the latter fact is remarkable:
 it is known,  that  the stable isomorphism group of a non-stationary  toric $AF$-algebra 
admits a faithful representation into the  matrix group $GL_n({\Bbb Z})$ \cite{E}.  
This elementary observation  implies   the Harvey conjecture.   
\begin{thm}\label{thm1}
For every surface $X$ of  genus $g\ge 2$,  there exists a faithful representation 
$\rho: Mod~(X)\rightarrow GL_{6g-6}({\Bbb Z})$.
\end{thm}
The structure of the article is as follows. 
The notation and facts used for the proof of theorem \ref{thm1} 
are introduced in Section 2.  Theorem \ref{thm1}  is  proved in Section 3.

\section{Preliminaries}
\subsection{$AF$-algebras}
{\bf A. The $C^*$-algebras.}
By a {\it $C^*$-algebra} one understands the Banach
algebra with an involution. Namely, a $C^*$-algebra
$A$ is an algebra over the complex numbers $\Bbb C$ with a norm $a\mapsto ||a||$
and an involution $a\mapsto a^*, a\in A$, such that $A$ is
complete with the respect to the norm, and such that
$||ab||\le ||a||~||b||$ and $||a^*a||=||a||^2$ for every
$a,b\in A$. If  $A$ is commutative, then the Gelfand
theorem says that $A$ is isometrically $*$-isomorphic
to the $C^*$-algebra $C_0(X)$ of the continuous complex-valued
functions on a locally compact Hausdorff space $X$.
For otherwise, the algebra $A$ represents a noncommutative topological
space.

\medskip\noindent
{\bf B.  The stable isomorphisms of $C^*$-algebras.} 
Let $A$ be a $C^*$-algebra deemed as a noncommutative
topological space. One can ask, when two such topological spaces
$A,A'$ are homeomorphic? To answer the question, let us recall
the topological $K$-theory.  If $X$ is a (commutative) topological
space, denote by $V_{\Bbb C}(X)$ an abelian monoid consisting
of the isomorphism classes of the complex vector bundles over $X$
endowed with the Whitney sum. The abelian monoid $V_{\Bbb C}(X)$
can be made to an abelian group, $K(X)$, using the Grothendieck
completion. The covariant functor $F: X\to K(X)$ is known to
map the homeomorphic topological spaces $X,X'$ to the isomorphic
abelian groups $K(X), K(X')$. Let  $A,A'$ be the $C^*$-algebras. If one wishes to
define a homeomorphism between the noncommutative topological spaces $A$ and $A'$, 
it will suffice to define an isomorphism between the abelian monoids $V_{\Bbb C}(A)$
and $V_{\Bbb C}(A')$ as suggested by the topological $K$-theory. 
The role of the complex vector bundle of the degree $n$ over the 
$C^*$-algebra $A$ is played by a $C^*$-algebra $M_n(A)=A\otimes M_n$,
i.e. the matrix algebra with the entries in $A$.  The abelian monoid
$V_{\Bbb C}(A)=\cup_{n=1}^{\infty} M_n(A)$ replaces the monoid 
$V_{\Bbb C}(X)$ of the topological $K$-theory. Therefore, 
the noncommutative topological spaces $A,A'$ are homeomorphic,
if the abelian monoids $V_{\Bbb C}(A)\cong V_{\Bbb C}(A')$ are isomorphic. 
The latter equivalence is called a {\it stable isomorphism}
of the $C^*$-algebras $A$ and $A'$ and is formally written as 
$A\otimes {\cal K}\cong A'\otimes {\cal K}$, where 
${\cal K}=\cup_{n=1}^{\infty}M_n$ is the $C^*$-algebra of compact
operators.  Roughly speaking, the stable isomorphism between
the $C^*$-algebras  means that they  are
homeomorphic as the noncommutative topological spaces.

\medskip\noindent
{\bf C. The $AF$-algebras.}
An {\it $AF$-algebra}  (approximately finite $C^*$-algebra) is defined to
be the  norm closure of an ascending sequence of the finite dimensional
$C^*$-algebras $M_n$'s, where  $M_n$ is the $C^*$-algebra of the $n\times n$ matrices
with the entries in ${\Bbb C}$. Here the index $n=(n_1,\dots,n_k)$ represents
a semi-simple matrix algebra $M_n=M_{n_1}\oplus\dots\oplus M_{n_k}$.
The ascending sequence mentioned above  can be written as 
$M_1\buildrel\rm\varphi_1\over\longrightarrow M_2
   \buildrel\rm\varphi_2\over\longrightarrow\dots,
$
where $M_i$ are the finite dimensional $C^*$-algebras and
$\varphi_i$ the homomorphisms between such algebras.  The set-theoretic limit
$A=\lim M_n$ has a natural algebraic structure given by the formula
$a_m+b_k\to a+b$; here $a_m\to a,b_k\to b$ for the
sequences $a_m\in M_m,b_k\in M_k$.  
The homomorphisms $\varphi_i$ can be arranged into  a graph as follows. 
Let  $M_i=M_{i_1}\oplus\dots\oplus M_{i_k}$ and 
$M_{i'}=M_{i_1'}\oplus\dots\oplus M_{i_k'}$ be 
the semi-simple $C^*$-algebras and $\varphi_i: M_i\to M_{i'}$ the  homomorphism. 
One has the two sets of vertices $V_{i_1},\dots, V_{i_k}$ and $V_{i_1'},\dots, V_{i_k'}$
joined by the $a_{rs}$ edges, whenever the summand $M_{i_r}$ contains $a_{rs}$
copies of the summand $M_{i_s'}$ under the embedding $\varphi_i$. 
As $i$ varies, one obtains an infinite graph called a {\it Bratteli diagram} of the
$AF$-algebra. The Bratteli diagram defines a unique  $AF$-algebra.

\medskip\noindent
{\bf D. The stationary $AF$-algebras.}
If the homomorphisms $\varphi_1 =\varphi_2=\dots=Const$ in the definition of 
the $AF$-algebra $A$,  the $AF$-algebra $A$ is called {\it stationary}. 
The Bratteli diagram of a stationary $AF$-algebra looks like a periodic 
graph  with the incidence matrix $A=(a_{rs})$ repeated over and over again. 
Since  matrix $A$ is a non-negative integer matrix, one can take a power of
$A$ to obtain a strictly positive integer matrix -- which we always assume 
to be the case.  The stationary $AF$-algebra has a non-trivial group
of the automorphisms \cite{E}, Ch.6.

\subsection{The Jacobi-Perron continued fraction}
{\bf A. The regular continued fractions.}
Let $a_1,a_2\in {\Bbb N}$ such that $a_2\le a_1$. Recall that the greatest common
divisor of $a_1$ and $a_2$, $GCD(a_1,a_2)$, can be determined from the Euclidean algorithm:
$$
\left\{
\begin{array}{cc}
a_1 &= a_2b_1 +r_3\nonumber\\
a_2 &= r_3b_2 +r_4\nonumber\\
r_3 &= r_4b_3 +r_5\nonumber\\
\vdots & \nonumber\\
r_{k-3} &= r_{k-2}b_{k-1}+r_{k-1}\nonumber\\
r_{k-2} &= r_{k-1}b_k,
\end{array}
\right.
$$
where $b_i\in {\Bbb N}$ and $GCD(a_1,a_2)=r_{k-1}$. 
The Euclidean algorithm can be written as the regular continued 
fraction
$$
\theta={a_1\over a_2}=b_1+{1\over\displaystyle b_2+
{\strut 1\over\displaystyle +\dots+ {1\over b_k}}}
=[b_1,\dots b_k].
$$
If $a_1$ and $a_2$ are non-commensurable, in the sense that $\theta\in {\Bbb R}-{\Bbb Q}$,
then the Euclidean algorithm never stops and $\theta=[b_1, b_2, \dots]$. Note that the regular  
continued fraction can be written in the matrix form:
$$
\left(\matrix{1\cr \theta}\right)=
\lim_{k\to\infty} \left(\matrix{0 & 1\cr 1 & b_1}\right)\dots
\left(\matrix{0 & 1\cr 1 & b_k}\right)
\left(\matrix{0\cr 1}\right). 
$$

\medskip\noindent
{\bf B. The Jacobi-Perron continued fractions.}
The Jacobi-Perron algorithm and connected (multidimensional) continued 
fraction generalizes the Euclidean algorithm to the case $GCD(a_1,\dots,a_n)$
when $n\ge 2$.  Specifically, let $\lambda=(\lambda_1,\dots,\lambda_n)$,
$\lambda_i\in {\Bbb R}-{\Bbb Q}$ and  $\theta_{i-1}={\lambda_i\over\lambda_1}$  with
$1\le i\le n$.   The continued fraction 
$$
\left(\matrix{1\cr \theta_1\cr\vdots\cr\theta_{n-1}} \right)=
\lim_{k\to\infty} 
\left(\matrix{0 &  0 & \dots & 0 & 1\cr
              1 &  0 & \dots & 0 & b_1^{(1)}\cr
              \vdots &\vdots & &\vdots &\vdots\cr
              0 &  0 & \dots & 1 & b_{n-1}^{(1)}}\right)
\dots 
\left(\matrix{0 &  0 & \dots & 0 & 1\cr
              1 &  0 & \dots & 0 & b_1^{(k)}\cr
              \vdots &\vdots & &\vdots &\vdots\cr
              0 &  0 & \dots & 1 & b_{n-1}^{(k)}}\right)
\left(\matrix{0\cr 0\cr\vdots\cr 1} \right),
$$
where $b_i^{(j)}\in {\Bbb N}\cup\{0\}$, is called the {\it Jacobi-Perron
algorithm (JPA)}. Unlike the regular continued fraction algorithm,
the JPA may diverge for certain vectors $\lambda\in {\Bbb R}^n$. However, 
for points of a generic subset of ${\Bbb R}^n$, the JPA converges. 
The convergence of the JPA algorithm can be characterized in terms of
the measured foliations. Let ${\cal F}\in\Phi_X$ be a measured foliation
on the surface $X$ of genus $g\ge 1$. Recall that ${\cal F}$ is called uniquely
ergodic if every invariant measure of ${\cal F}$ is a multiple
of the Lebesgue measure. By the Masur-Veech theorem, there exists
a generic subset $V\subset \Phi_X$ such that each ${\cal F}\in V$
is uniquely ergodic \cite{Mas1}, \cite{Vee1}.
We let $\lambda=(\lambda_1,\dots,\lambda_{n})$ be the vector with
coordinates $\lambda_i=\mu ({\gamma_i})$, where $\gamma_i\in H_1^{odd}(\widetilde S)$;
by an abuse of notation, we shall say that $\lambda\in V$. 
In view of  duality  between the measured foliations and the interval exchange 
transformations \cite{Mas1}, the JPA converges if and only if 
$\lambda\in V\subset {\Bbb R}^{n}$ \cite{Bau1}.

\section{Proof of theorem 1}
As before, let $W$ denote  the set of toric $AF$-algebras of genus $g\ge 2$. 
Let $G$ be a finitely generated group and  $G\times W\to W$ be an action of $G$ on $W$ 
by  the stable isomorphisms of toric $AF$-algebras; in other words, 
 $\gamma ({\Bbb A}_{\theta})\otimes {\cal K}\cong
{\Bbb A}_{\theta}\otimes {\cal K}$  for all $\gamma\in G$ and all 
${\Bbb A}_{\theta}\in W$.   The following preparatory lemma will be 
important.
\begin{lem}\label{lm2}
For each ${\Bbb A}_{\theta}\in W$,  there exists a 
representation $\rho_{{\Bbb A}_{\theta}}: G\to GL_{6g-6}({\Bbb Z})$.  
\end{lem}
{\it Proof.}
The proof of  lemma is based on the following well known criterion
of the stable isomorphism for the $AF$-algebras: a pair of such algebras
${\Bbb A}_{\theta}, {\Bbb A}_{\theta'}$ are stably isomorphic if and only
if their Bratteli diagrams coincide, except (possibly) a finite part
of the diagram, see \cite{E}, Theorem 2.3. (Note, that the order isomorphism
between  the dimension groups mentioned in the original text, can be reformulated 
in the language of the Bratteli diagrams as stated.)

Let $G$ be a finitely generated
group on the generators  $\{\gamma_1, \dots, \gamma_m\}$ and ${\Bbb A}_{\theta}\in W$.
Since $G$ acts on the toric $AF$-algebra ${\Bbb A}_{\theta}$ by stable isomorphisms,
the toric $AF$-algebras ${\Bbb A}_{\theta_1}:=\gamma_1({\Bbb A}_{\theta}),\dots,   
{\Bbb A}_{\theta_m}:=\gamma_m({\Bbb A}_{\theta})$ are stably isomorphic to 
${\Bbb A}_{\theta}$; moreover, by transitivity, they are also pairwise stably isomorphic.
Therefore, the Bratteli diagrams of ${\Bbb A}_{\theta_1},\dots, {\Bbb A}_{\theta_m}$ 
 coincide everywhere except, possibly, some finite parts.  
We shall denote by ${\Bbb A}_{\theta_{\max}}\in W$
a toric $AF$-algebra, whose Bratteli diagram is the maximal common part
of the Bratteli diagrams of ${\Bbb A}_{\theta_i}$ for $1\le i\le m$;
such a choice is unique and defined correctly because the set $\{{\Bbb A}_{\theta_i}\}$
is a finite set.  By the definition of a toric $AF$-algebra, the vectors 
$\theta_i=(1,\theta_1^{(i)},\dots,\theta_{6g-7}^{(i)})$
are related to the vector $\theta_{\max}=(1, \theta_1^{(\max)},\dots,\theta_{6g-7}^{(\max)})$ 
by the formula:
$$
\left(\matrix{1\cr \theta_1^{(i)}\cr\vdots\cr\theta_{6g-7}^{(i)}} \right)
=\underbrace{
\left(\matrix{0 &  0 & \dots & 0 & 1\cr
              1 &  0 & \dots & 0 & b_1^{(1)(i)}\cr
              \vdots &\vdots & &\vdots &\vdots\cr
              0 &  0 & \dots & 1 & b_{6g-7}^{(1)(i)}}\right)
\dots 
\left(\matrix{0 &  0 & \dots & 0 & 1\cr
              1 &  0 & \dots & 0 & b_1^{(k)(i)}\cr
              \vdots &\vdots & &\vdots &\vdots\cr
              0 &  0 & \dots & 1 & b_{6g-7}^{(k)(i)}}\right)
}_{A_i}
\left(\matrix{1\cr \theta^{(\max)}_1\cr\vdots\cr\theta^{(\max)}_{6g-7}} \right)
$$
The above expression can be written in the matrix form $\theta_i=A_i\theta_{\max}$, where 
$A_i\in GL_{6g-6}({\Bbb Z})$. Thus, one gets a matrix representation of the
generator $\gamma_i$,  given by the formula $\rho_{{\Bbb A}_{\theta}}(\gamma_i):=A_i$.
The map $\rho_{{\Bbb A}_{\theta}}: G\to GL_{6g-6}({\Bbb Z})$ extends to the rest of the group $G$
 via its values on the generators;  namely,  for every $g\in G$ one sets $\rho_{{\Bbb A}_{\theta}}(g)= A_1^{k_1}\dots A_m^{k_m}$,
whenever $g=\gamma_1^{k_1}\dots \gamma_m^{k_m}$. It is verified (by induction), that  the map
$\rho_{{\Bbb A}_{\theta}}: G\to GL_{6g-6}({\Bbb Z})$ is a homomorphism, since 
$\rho_{{\Bbb A}_{\theta}}(g_1g_2)=\rho_{{\Bbb A}_{\theta}}(g_1)\rho_{{\Bbb A}_{\theta}}(g_2)$
for $\forall g_1,g_2\in G$. 
  Lemma \ref{lm2} follows.
$\square$

\medskip
Let  $W_{aper}\subset W$ be a set consisting  of the toric $AF$-algebras, whose  Bratteli diagrams
do not contain periodic (infinitely repeated) blocks;  these are known as  non-stationary toric  $AF$-algebras
and they are generic in the set $W$ endowed with the natural topology. 
We call  the action of $G$  on the toric $AF$-algebra  ${\Bbb A}_{\theta}\in W$
{\it free}, if   $\gamma ({\Bbb A}_{\theta})={\Bbb A}_{\theta}$ implies $\gamma=Id$.
\begin{lem}\label{lm3}
Let ${\Bbb A}_{\theta}\in W_{aper}$  and  $G$ be  
free on the ${\Bbb A}_{\theta}$.  Then  $\rho_{{\Bbb A}_{\theta}}$ is a 
faithful representation.
\end{lem}
{\it Proof.}
Since the action of $G$ is free, to prove that  $\rho_{{\Bbb A}_{\theta}}$  is faithful, 
it remains  to show, that in the formula  $\theta_i=A_i\theta_{\max}$, it holds
$A_i=I$, if and only if,  $\theta_i=\theta_{\max}$, where $I$
is the unit matrix.   Indeed,  it is immediate that $A_i=I$ implies $\theta_i=\theta_{\max}$.
Suppose now that  $\theta_i=\theta_{\max}$ and, let to the contrary, $A_i\ne I$. 
One gets $\theta_i=A_i \theta_{\max}=\theta_{\max}$.  Such an equation has  a non-trivial solution, 
if and only if,  the vector $\theta_{\max}$ has a periodic 
Jacobi-Perron fraction;  the period of  such a fraction  is given by the matrix $A_i$. This 
is impossible, since it has been assumed, that ${\Bbb A}_{\theta_{\max}}\in W_{aper}$.
The contradiction finishes the proof of lemma \ref{lm3}.
$\square$

\medskip
Let $G=Mod~(X)$, where $X$ is a surface of genus $g\ge 2$. The group $G$
is finitely generated \cite{Deh1}; it  acts on the Teichm\"uller space $T(g)$
by isomorphisms of the Riemann surfaces. Moreover, the action of $G$ is free on a 
generic set, $U\subset T(g)$,  consisting of the Riemann surfaces with the trivial group 
of the automorphisms.  Let $F: V\to W$ be the Teichm\"uller functor  between the Riemann surfaces and 
toric $AF$-algebras (lemma \ref{lm0});  the following  is true.
\begin{lem}\label{lm4}
The pre-image  $F^{-1}(W_{aper})$ is a generic set in the space  $T(g)$.
\end{lem}
{\it Proof.}
Note,  that the set of stationary toric $AF$-algebras is a countable set.
The functor $F$ is a surjective map,  which is continuous with respect to the natural topology 
on the sets $V$ and $W$.  Therefore,  the pre-image of the complement of a countable set is a generic set. 
$\square$

\medskip
To finish the proof,  consider the set $U\cap F^{-1}(W_{aper})$;
this set is a non-empty set, since it is the intersection of  two generic subsets of $T(g)$.
Let $R$ be a point (a Riemann surface) in the above set.  In  view of lemma \ref{lm0},
group $G$ acts on the toric $AF$-algebra ${\Bbb A}_{\theta}=F(R)$ by the stable
isomorphisms.  By the construction, the action is free and ${\Bbb A}_{\theta}\in W_{aper}$.
In view of lemma \ref{lm3},  one gets a faithful representation $\rho=\rho_{{\Bbb A}_{\theta}}$
of the group $G=Mod~(X)$ into the matrix group $GL_{6g-6}({\Bbb Z})$.
Theorem \ref{thm1} follows.
$\square$



\vskip1cm

\textsc{The Fields Institute for Mathematical Sciences, Toronto, ON, Canada,  
E-mail:} {\sf igor.v.nikolaev@gmail.com}

\smallskip
{\it Current address: 101-315 Holmwood Ave., Ottawa, ON, Canada, K1S 2R2}  

\end{document}